\documentclass[a4paper,10pt]{article} \usepackage{amsmath,amssymb,amsthm}
\usepackage[english]{babel}

\newtheorem{theo}{Theorem}

\newtheorem{lemm}[theo]{Lemma} 
\newtheorem{prop}[theo]{Proposition}

\newcommand{\Ra}{\mathbb R}                   

\newcommand{\scal}[1]{\langle #1 \rangle}

\newcommand{\finpreuve}{\hfill $\Box$}

\newcommand{\name}{$\underline{\qquad \qquad}$}

\begin{document}

\author{  Jean-Marc Bouclet}

\title{Normal form of the metric for a class of Riemannian manifolds with ends}

\maketitle

\begin{abstract} 
In many problems of PDE involving the Laplace-Beltrami operator on manifolds with ends, it is often useful to introduce radial or geodesic normal coordinates near infinity.
In this paper, we prove the existence of such coordinates for a general class of manifolds with ends, which includes asymptotically conical and hyperbolic manifolds. We  study the decay rate to the metric at infinity associated to radial coordinates and also show that the latter metric is always conformally equivalent to the metric at infinity associated to the original coordinate system. We finally give several examples illustrating the sharpness of our results.
\end{abstract}

\noindent {\bf Keywords:} Manifolds with ends, radial coordinates, geodesic normal coordinates.

\section{Introduction and result}
\setcounter{equation}{0}

The purpose of this note is to study the existence and some properties of radial (or geodesic normal) coordinates at infinity  on manifolds with ends, for a general class of ends.
Our motivation comes from geometric spectral and scattering theory (see e.g. \cite{Melrose} for important aspects of this topic), but our results may be of independent interest. The kind of manifolds we consider is as follows. We assume that, away from a compact set, they are a finite union of ends $ {\mathcal E} $  isometric to $ \big( (R, +\infty) \times {\mathcal S} , {\mathbf G} \big) $ with ${\mathcal S}$ a compact manifold (of dimension $ n-1 \geq 1 $ in the sequel) and ${\mathbf G}$ of the form
\begin{eqnarray} 
 {\mathbf G} = {\mathbf a} dx^2 +  2 {\mathbf b}_i   dx  d \theta_i / w (x) + {\mathbf g}_{ij} d \theta_i d \theta_j / w (x)^2 , \label{metriquedirecte}
\end{eqnarray} 
(using the summation convention) with coefficients satisfying, as $ x \rightarrow \infty $,
\begin{eqnarray}
 {\mathbf a} (x,\theta) \rightarrow 1 , \qquad {\mathbf b}_i (x,\theta) \rightarrow 0, \qquad {\mathbf g}_{ij}(x,\theta) \rightarrow \overline{\mathbf g}_{ij}(\theta)=: \overline{\mathbf g} \left( \frac{\partial}{\partial \theta_i} , \frac{\partial}{\partial \theta_j} \right) . \label{nonquantitatif}
\end{eqnarray}
The nature of the end is determined by the function $ w $ which we assume here to be positive, smooth and, more importantly, 
$$ w(x) \rightarrow 0 \qquad x \rightarrow + \infty , $$ meaning that we consider large ends. The two main important examples are asymptotically conical manifolds (or scattering manifolds) for which $ w (x) = x^{-1} $ and asymptotically hyperbolic manifolds for which $ w (x) = e^{-c x} $ for some $c>0$. 
In (\ref{nonquantitatif}), $ \theta_{\mathcal S} = \big( \theta_1 , \ldots , \theta_{n-1} \big) : U \subset {\mathcal S} \rightarrow \Ra^{n-1} $ are local coordinates  on ${\mathcal S}$ so if  $ \pi : {\mathcal E} \rightarrow {\mathcal S} $ is the projection,  we obtain local coordinates on $ {\mathcal E} $ by considering $ (x, \theta_1 \circ \pi, \ldots , \theta_{n-1} \circ \pi) $ which, for simplicity of the notation, we denote by  $ (x,\theta_1 , \ldots , \theta_{n-1} ) $. The precise meaning of (\ref{nonquantitatif}) is that the convergence holds in $ C^{\infty} \big( \theta_{\mathcal S} (U) \big) $; such a statement is intrinsic in that it is invariant under the change of coordinates on $ {\mathcal S} $. We call $ \overline{\mathbf g} $ the metric at infinity with respect to this product decomposition.

For analytical purposes, it is often very useful to work in a system of coordinates such that
$ {\mathbf a} \equiv 1 $ and $ {\mathbf b}_i \equiv 0 $, i. e. to replace $x$ by a new coordinate $ t $ such that
\begin{eqnarray}
 {\mathbf G} = dt^2 + {\mathbf h}_{ij} d \theta_i d \theta_j / w (t)^2, \qquad {\mathbf h}_{ij}(t,\theta) \rightarrow \overline{\mathbf h}_{ij}(\theta)=: \overline{\mathbf h} \left( \frac{\partial}{\partial \theta_i} , \frac{\partial}{\partial \theta_j} \right) \ \ \mbox{as} \ t \rightarrow + \infty , \label{normalformintro} 
\end{eqnarray}
at the expense of changing $ \overline{\mathbf g} $ into a possibly different metric $  \overline{\mathbf h} $. One then says that $t$ is a radial coordinate (see for instance \cite{Kumura} for the terminology). Using such coordinates, the Laplacian can then be reduced, up to conjugation by a suitable function, to an operator of the form $ - \partial_t^2 + Q (t) $ with $ Q (t) $ an elliptic operator on $ {\mathcal S}$  asymptotic to $ - w(t)^2 \Delta_{\overline{\mathbf h}} $ as $t \rightarrow \infty $ (see e.g. (1.1) in \cite{CaVo2}). The absence of crossed term of the form $ \partial_t \partial_{\theta_i} $ is convenient for Born-Oppenheimer approaches, i. e. to consider $ - \partial_t^2 + Q (t) $ as a one dimensional Schr\"odinger operator with an operator valued potential (see for instance \cite{BouVP} for applications in this spirit); in the special case when $ Q (t) $ is exactly $  - w (t)^2 \Delta_{{\overline{\mathbf h}}} $, i. e. if $ {\mathbf G} = dt^2 + \overline{\mathbf h} / w (t)^2 $, one can use separation of variables as is well known. Important questions requiring such a reduction of the metric also include resolvent estimates \cite{Burq,CaVo1,CaVo2} (construction of Carleman weights) or inverse problems \cite{JoshiSaBaretto1,JoshiSaBaretto2} (reduction to a problem on $ {\mathcal S} $). 

In the works \cite{Burq,CaVo1,CaVo2,JoshiSaBaretto1,JoshiSaBaretto2}, the reduction of $ {\mathbf G} $ to the normal form (\ref{normalformintro}) is either proved on particular cases \cite{Burq,JoshiSaBaretto1} (conical ends) and \cite{JoshiSaBaretto2} (asymptotically hyperbolic ends), or even taken as an assumption in \cite{CaVo1,CaVo2}. For this reason and also in the perspective of studying intermediate models between the conical and the asymptotically hyperbolic cases, we feel  worth proving in detail the existence of radial coordinates for general manifolds with ends (i. e. associated to $w$ satisfying the assumption (\ref{hypothesew}) below). Another motivation is that, although the existence of radial coordinates may seem intuitively  clear, there are some subtleties on the rate of convergence to the asymptotic metric. We shall in particular show that, even if the convergences in (\ref{nonquantitatif}) are fast as $x \rightarrow \infty $, it may happen that the decay in radial coordinates, i. e. the rate of convergence to $ \overline{\mathbf h} $ in (\ref{normalformintro}), is slow. We shall see how this depends on  $w$. This point is important in scattering theory since it means  that the reduction to (\ref{normalformintro}) may be at the price of considering a long range type of decay. As a last point, we shall also describe the relationship between $ \overline{\mathbf g} $ and $ \overline{\mathbf h} $. For the class of functions $w$ we are going to consider, we shall see that $ \overline{\mathbf h} $ is always conformally equivalent to $ \overline{\mathbf g} $, as is well known in the asymptotically hyperbolic case. In certain situations, such as the conical case, the conformal factor is equal to 1 (i. e. there is no conformal change) and this will be covered by our result.

Let us now state our main result precisely.

 First, for simplicity and without loss of generality, we will assume that $ {\mathcal M} = {\mathcal E} = (R,\infty) \times {\mathcal S} $ equipped with a Riemannian metric $ {\mathbf G} $ as in (\ref{metriquedirecte}). We will use a quantitative version of (\ref{nonquantitatif}) given in term of symbol classes $ S^m $.  Recall that, given $ m \in \Ra$ and a function $f$ defined  on a semi-infinite interval $ (M,+\infty) $ or on $ (M,+\infty) \times V $, with $V$ an open subset of $ \Ra^{n-1} $, we have
 $$ f \in S^m  \qquad \stackrel{\rm def}{ \Longleftrightarrow} \qquad \partial_x^j \partial_{\theta}^{\alpha} f = {\mathcal O} \big( \scal{x}^{m - j} \big) , $$
 on $ (M,+ \infty) \times K $ for all $ K \Subset V $. Occasionally we shall also say that a function or a tensor defined on $ (M,+\infty) \times {\mathcal S} $ belongs to  $ S^m $ if its pullback by every coordinate chart of an atlas of $ {\mathcal S} $ is in $ S^m $.  
 
The precise assumptions on $ {\mathbf G} $  are as follows.
We assume first that, for some $ \lambda > 0 $ and $ \varepsilon > 0 $, 
\begin{eqnarray}
w \in S^{-\lambda} , \qquad  \left( \frac{w^{\prime}}{w} \right)^{\prime} \in S^{-1- \varepsilon} , \label{hypothesew}
\end{eqnarray}
where $ S^{m} = S^m (R,\infty) $ for $ m = - \lambda $ and $ - 1 - \varepsilon $.
The condition on $ (w^{\prime}/w)^{\prime} $ implies the existence of the non positive real number
\begin{eqnarray}
 \kappa : = \lim_{x \rightarrow + \infty} \frac{w^{\prime}(x)}{w(x)} . \label{limitew}
\end{eqnarray}
Notice that $ \kappa \leq 0 $. Otherwise $ w^{\prime} $ should be positive at infinity hence $w$ should be increasing which would be incompatible with the fact that $ w \in S^{-\lambda} $ (recall that $ w > 0 $).
To state our second assumption,  we set $ {\mathbf b} = ({\mathbf b}_1 , \ldots , {\mathbf b}_{n-1}) $, $ {\mathbf g} = ({\mathbf g}_{ij}) $ and $ \overline{\mathbf g} = (\overline{\mathbf g}_{ij}) $ (see (\ref{metriquedirecte}) and (\ref{nonquantitatif})). We assume that
\begin{eqnarray}
 {\mathbf a} - 1 \in S^{-\mu}  , \qquad {\mathbf b} \in S^{-\nu} , \qquad  {\mathbf g} - \overline{\mathbf g} \in S^{-\tau}  , \label{quantitatif}
\end{eqnarray}
where  $ S^{m} = S^m ((R,\infty) \times \theta_{\mathcal S}(U)) $ (for all charts $ \theta_{\mathcal S} : U \rightarrow \theta_{\mathcal S}(U) $ of some atlas of $ {\mathcal S} $) and with exponents satisfying
\begin{eqnarray}
 \mu \geq 1 + \tau , \qquad  \nu  \geq \frac{1 + \tau}{2} , \qquad  \lambda \geq  \frac{1 + \tau}{2} , \qquad \mbox{with} \ \tau > 0. \label{conditionexposant}
\end{eqnarray}
 We finally  define the outgoing normal geodesic flow. Given $ r > R $, denote by $ \nu_r $ the outgoing normal vector field to the hypersurface $ \{ r \} \times {\mathcal S} \subset {\mathcal M} $. Here outgoing means that $ \langle d x , \nu_r
 \rangle > 0 $. The outgoing normal geodesic flow $ N_r $ is then 
 $$ N_r (t,\omega)  := \exp_{(r,\omega)} (t \nu_r) , \qquad \omega \in {\mathcal S} , \ t \geq 0 , $$
 namely the exponential map on $ {\mathcal M} $ with starting point on $ \{ r\} \times {\mathcal S} $, initial speed $ \nu_r $ and nonnegative time.
 \begin{theo} \label{theoremecoordonnees} Assume (\ref{hypothesew}), (\ref{quantitatif}) and (\ref{conditionexposant}). Then, for all $ r \gg 1 $, $N_r$ has the following properties.
 \begin{enumerate} \item{It is complete in the future (i. e. is defined for all $ t \geq 0 $).}
 \item{It is a homeomorphism (resp. a diffeomorphism) between $ [0,\infty )_t \times {\mathcal S} $ (resp. $ (0,\infty) \times {\mathcal S} $) and $ [r,\infty )_x \times {\mathcal S} $ (resp. $ (r,\infty) \times {\mathcal S} $). }
 \item{There exists a diffeomorphism $ \Omega_r : {\mathcal S} \rightarrow {\mathcal S} $ and a real function $ \phi_r : {\mathcal S} \rightarrow \Ra $ such that
 $$ N_r^* {\mathbf G} = d t^2 + w (t)^{-2} {\mathbf h} (t) $$ with $ \big( {\mathbf h} (t) \big)_{t > 0} $ a smooth family of metrics on $ {\mathcal S} $ such that 
 \begin{eqnarray}
  {\mathbf h} (t) -  \overline{\mathbf h} \in S^{-\min (\tau , \varepsilon)} , \qquad \mbox{with} \qquad  \overline{\mathbf h} := e^{- 2 \kappa \phi_r} \Omega_r^* \overline{\mathbf g} . \label{formuleaveckappa}
 \end{eqnarray}}
 \end{enumerate}
 \end{theo}
 
 Note the dependence on $ \kappa $ in (\ref{formuleaveckappa}). In particular, if $ \kappa = 0 $, there is no conformal factor. 
Observe also that the decay rate  of $ {\mathbf h} - \overline{\mathbf h} $ in (\ref{formuleaveckappa}) can in principle be worse than the one of  $ {\mathbf g} - \overline{\mathbf g} $ in (\ref{quantitatif}). We shall see that this can be the case in some of the examples below.

\bigskip

\noindent {\bf Examples.} {\it 1. Asymptotically conical metrics:} $ w (x) = x^{-1} $ (for $ x > R > 0$).  We have obviously
$$ \lambda = 1, \qquad \varepsilon = 1, \qquad \kappa = 0 . $$
On one hand $ \kappa = 0 $, so the metric at infinity is not affected by a conformal factor, but on the other hand $ \varepsilon = 1 $ so $ {\mathbf h} (t) $ is in general a long range perturbation of $ \overline{\mathbf h} $.  Actually, one can see that
\begin{eqnarray}
 {\mathbf h}(t) = (1 + 2 \phi_r t^{-1})  \overline{\mathbf h} + o (t^{-1}) , \label{expansionexample1}
\end{eqnarray} 
which shows that the decay rate of $ {\mathbf h}(t) -  \overline{\mathbf h} $ is only $ S^{-1} $ (see the proof of Theorem \ref{theoremecoordonnees} below in Subsection \ref{subsection2point2} for a justification of (\ref{expansionexample1})).
\medskip

\noindent {\it 2. Asymptotically hyperbolic metrics:} $ w (x) = e^{- c x} $ (with $ c > 0 $). In this case, we can take
$$ \lambda > 0 \ \mbox{ arbitrarily large}, \qquad \varepsilon > 0 \ \mbox{ arbitrarily large}, \qquad \kappa = c . $$
Here $ \kappa \ne 0 $ hence the metric at infinity is only conformally equivalent to the original one. On the other hand, since $ \varepsilon $ can be taken as large as we wish, in particular larger than $ \tau $, the decay rate of $ {\mathbf h} $ to $ \overline{\mathbf h} $ cannot be worse than the one of $  {\mathbf g} $ to $ \overline{\mathbf g} $ in (\ref{quantitatif}). 
\medskip

\noindent {\it 3. An intermediate case.} For the function $w (x) = e^{-x-x^{\beta}}$, with $ 0 < \beta < 1 $, we have
$$ \lambda > 0 \ \mbox{ arbitrarily large}, \qquad \varepsilon = 1 - \beta, \qquad \kappa = 1 . $$  
This suggests  that both a conformal factor and a weaker decay (if $ \varepsilon < \tau $)  happen at the same time. Actually the decay can indeed be weaker if $ \varepsilon < \tau $,  for one can show that
\begin{eqnarray}
{\mathbf h} (t) = (1 + 2 \beta \phi_r t^{\beta-1}) \overline{\mathbf h} + o (t^{\beta-1}) + O (t^{-\tau})  .\label{expansionexample3}
\end{eqnarray}
See again  the proof of Theorem \ref{theoremecoordonnees} below for a justification of this expansion.
\section{The outgoing normal geodesic flow}
\setcounter{equation}{0}

\subsection{The main estimates}
In this subsection, we fix some notation and state intermediate results leading fairly directly to Theorem \ref{theoremecoordonnees} which is proved in Subsection \ref{subsection2point2}. The more technical proofs are postponed to the next sections.

It will be convenient to use some fixed geodesic distance $ d(.,.) $ on $ {\mathcal S} $ associated to an arbitrary Riemannian metric (which has nothing to do with $ \overline{\mathbf g} $). We then fix a cover of ${\mathcal S}$  by finitely many coordinates patches. At any $ \omega_0 \in {\mathcal S} $, there is a chart $ \theta_{\mathcal S} : U \subset {\mathcal S} \rightarrow V \subset \Ra^{n-1} $ and, if we set $ \theta_0 = \theta_{\mathcal S} (\omega_0) $, there is $ \epsilon_{\omega_0} $ such that 
\begin{eqnarray}
 B (\theta_0 , 4 \epsilon_{\omega_0}) \Subset V , \label{conditionreconnaissance}
\end{eqnarray} 
where, here and below, the ball $ B (\theta_0, \epsilon) $ refers to a fixed norm $ | \cdot | $ on $ \Ra^{n-1} $.  By the compactness of ${\mathcal S}$, we have
\begin{eqnarray}
 {\mathcal S} = \bigcup_{\omega_0 \in {\rm finite \ set}} \theta_{\mathcal S}^{-1} \big( B (\theta_0, \epsilon_{\omega_0}) \big) . \label{recouvrement}
\end{eqnarray}
Furthermore, we can assume that, for some fixed $ C > 0 $ depending on $d$ and the cover (\ref{recouvrement}), 
\begin{eqnarray}
d (\omega,\omega^{\prime}) \leq C |\theta_{\mathcal S}(\omega) - \theta_{\mathcal S}(\omega^{\prime})|, \qquad \omega, \ \omega^{\prime} \in \theta_{\mathcal S}^{-1} \big( B (\theta_0, 3 \epsilon_{\omega_0}) \big) , \label{globalisationdistance2}
\end{eqnarray}
with $d$ the distance which was chosen above.

We next summarize the expressions of several important objects in the coordinate patch of $ {\mathcal M} $ associated to the patch $ \theta_{\mathcal S}^{-1}(B(\theta_0,4 \epsilon_{\omega_0})) $ of $ {\mathcal S} $.  We will study the geodesic flow through its hamiltonian expression on the cotangent bundle and thus need to compute the dual metric. To this end, we  recall that (\ref{metriquedirecte}) can be recast in matrix form as
\begin{eqnarray}
{\mathbf G} \equiv  \left( \begin{matrix}
1 & 0     \\
0 & w (x)   
\end{matrix}
 \right)^{-1} \left( \begin{matrix}
{\mathbf a} & {\mathbf b}^T    \\
{\mathbf b} & {\mathbf g}   
\end{matrix}
 \right) \left( \begin{matrix}
1 & 0     \\
0 & w (x) \end{matrix}
 \right)^{-1}. \label{matrixmetric}   
\end{eqnarray}
Then the dual metric, obtained by inverting (\ref{matrixmetric}),  is given by $  \left( \begin{matrix}
1 & 0     \\
0 & w (x)   
\end{matrix}
 \right) \left( \begin{matrix}
a &  b^T    \\
b & g   
\end{matrix}
 \right) \left( \begin{matrix}
1 & 0     \\
0 & w (x) \end{matrix}
 \right) $ with
 \begin{eqnarray}
a  =  \frac{1}{{\mathbf a} - {\mathbf b}^T {\mathbf g}^{-1} {\mathbf b}}, \qquad \qquad b  =  - a {\mathbf g}^{-1} {\mathbf b}, \qquad \qquad g  =  {\mathbf g}^{-1} + a {\mathbf g}^{-1} {\mathbf b} {\mathbf b}^T {\mathbf g}^{-1} .  \label{formule}
\end{eqnarray} 
Note that, by possibly increasing $R$ and  by (\ref{nonquantitatif}), we may assume that $ {\mathbf a} - {\mathbf b}^T {\mathbf g}^{-1} {\mathbf b} $ does not vanish. 
It is important to note that, by (\ref{quantitatif}) and (\ref{conditionexposant}), we have
\begin{eqnarray}
 a - 1 \in S^{- \min (\mu , 2 \nu)} \subset S^{-1-\tau}  , \qquad  b \in S^{-\nu} , \qquad  g - \bar{\mathbf g}^{-1} \in S^{-\tau} . \label{decroissancemetriquedualenew}
\end{eqnarray}
According to the notation (\ref{formule}), the dual metric, i. e. the principal symbol of the Laplacian, reads
\begin{eqnarray}
 p (x,\theta,\rho,\eta) := a (x,\theta)\rho^2 + 2 w(x) \rho b (x,\theta) \cdot \eta + w(x)^2 \eta \cdot g (x,\theta) \eta  , \label{formesymboletotal}
\end{eqnarray}
with $ \rho \in \Ra $ and $ \eta \in \Ra^{n-1} $. We  denote by $ \big( x^t,\theta^t, \rho^t , \eta^t \big) $ the hamiltonian flow  of $p$, namely the solution to
\begin{eqnarray}
\dot{x}^t = \frac{\partial p}{\partial \rho} , \qquad \dot{\theta}^t = \frac{\partial p}{\partial \eta}, \qquad  \dot{\rho}^t = - \frac{\partial p}{\partial x}
 \qquad \dot{\eta}^t = - \frac{\partial p}{\partial \theta}, \label{hamiltonequation}
\end{eqnarray}
with initial condition at $t=0$ to be specified. A simple calculation shows that the outgoing normal to  $ \{ r \} \times {\mathcal S} $ is the vector field
$$ \nu_r  =    a^{1/2} \frac{\partial}{\partial x} +   w (x) \frac{b}{ a^{1/2}} \cdot \frac{\partial}{\partial \theta}  , $$
where $ a $ and $ b $ are evaluated at $ (r,\theta) = (r,\theta_{\mathcal S}(\omega)) $.  The associated co-normal form $ \nu_r^{\flat} $, i. e. such that $ {\mathbf G} (\nu_r,.) = \nu_r^{\flat} $,  is then
$$ \nu_r^{\flat} = a^{-1/2}dx , $$
so the geodesic starting at $ (r, \omega ) $ with $ \nu_r  $ as initial velocity, i. e. $ \exp_{(r,\omega)}(t \nu_r) $, is given in these local coordinates by
\begin{eqnarray}
x_{\nu_r} (t, \theta ) : =  x^{t/2} \big(r , \theta , a^{-1/2}  , 0 \big), \qquad \theta_{\nu_r} (t,\theta) := \theta^{t/2} \big(r , \theta , a^{-1/2}  , 0 \big)    . 
\label{expressiongeodesique}
\end{eqnarray}
Here the factor $ 1/2 $ on the time is due to the fact that we consider the Hamiltonian flow of $ p $ rather than the one of $ p^{1/2} $. 
 We also note in passing  that the condition $ {\mathbf G} (\nu_r,\nu_r)=1 $ reads 
\begin{eqnarray}
p \big( r,\theta,a^{-1/2} ,0 \big) = 1 . \label{energieautomatique}
\end{eqnarray}
The expression of the normal geodesic flow given by (\ref{expressiongeodesique}) is of course meaningful only  as long as the geodesic remains in the coordinate patch. We  shall see below that, if $r$ is large enough and $  \theta $ is restricted to $ B (\theta_0 , 2 \epsilon_{\omega_0}) $ (which is technically more convenient than $ B (\theta_0, \epsilon_{\omega_0})$, though the latter would be sufficient by (\ref{recouvrement})), then the  geodesic remains in the same coordinate patch for all $t \geq 0$ (thus is complete in the future) and satisfies suitable estimates. To make the proof as clear as possible, we pick up its main steps in the following propositions which will be proved in separate subsections.

\begin{prop}[the geodesic flow in a chart] \label{propositiondifferee} Assume (\ref{hypothesew}), (\ref{quantitatif}) and (\ref{conditionexposant}).
Then, for all $ M > 1 $, there exists $  X >0 $ such that, for all initial condition $ (x,\theta,\rho,\eta) $ of (\ref{hamiltonequation}) satisfying
\begin{eqnarray}
 x \geq X, \qquad \theta \in B(\theta_0, 2 \epsilon_{\omega_0}), \qquad \rho \in \big[ M^{-1}, M \big], \qquad |\eta| \leq M ,  \label{conditionsinitiales2}
\end{eqnarray}
the hamiltonian flow of $p$ is defined for all $  t \geq 0 $ and satisfies
\begin{eqnarray}
x^t  \geq  x + \frac{t}{M},  \qquad \theta^t \in B (\theta_0 , 3 \epsilon_{\omega_0}) . \label{stabilitedomaine} 
\end{eqnarray}
Furthermore, for all $j\geq 1$ and all $ \partial^{\gamma} = \partial_x^k \partial_{\theta}^{\alpha} \partial_{\rho}^l \partial_{\eta}^{\beta} $, we have the estimates
\begin{eqnarray}
\big|\partial_t^{j} \partial^{\gamma} (x^t - x - 2 t p^{1/2} ) \big| & \lesssim & \scal{x+t}^{-\tau-j} , \label{flotpratique1} \\
\big| \partial_t^{j} \partial^{\gamma} (\theta^t - \theta)\big| & \lesssim & \scal{x+t}^{-\tau-j} \label{flotpratique2}
\end{eqnarray}
where $ p = p (x,\theta,\rho,\eta) $. 
\end{prop}

\noindent {\it Proof.} See Section \ref{subsectionproof}. 

\bigskip

We now derive here a proposition on the outgoing normal geodesic flow from which  Theorem \ref{theoremecoordonnees} will follow easily.  We introduce the notation
\begin{eqnarray} 
 N_r =: (x_r, \omega_r) \label{notationflotnormal}
\end{eqnarray}
for the components of $ N_r $ on $ (R,+\infty) $ and $ {\mathcal S} $, respectively. 
Note the relationship between (\ref{notationflotnormal}) and (\ref{expressiongeodesique}):
\begin{eqnarray}
  x_r   \left( t , \theta_{\mathcal S}^{-1} ( \theta ) \right) = x_{\nu_r} \big( t ,  \theta \big), \qquad \left( \theta_{\mathcal S} \circ \omega_r \right) \left( t , \theta_{\mathcal S}^{-1} ( \theta ) \right) =  \theta_{\nu_r}  \big( t , \theta \big) . \label{traductioncarte}
\end{eqnarray}

\begin{prop}[Global properties of the normal flow] \label{propnormalflow} For all $ r \gg 1 $, the following properties hold.
\begin{enumerate}
\item{For each $ t \geq 0 $, $ \omega_r (t,.) $ is a diffeomorphism from $ {\mathcal S} $ to $ {\mathcal S} $ and
$$ d \big( \omega , \omega_r (t, \omega) \big) \leq C \scal{r}^{-\tau}, \qquad r \gg 1, \ t \geq 0 , \ \omega \in {\mathcal S}  ,$$
with $ C $ independent of $ r,t,\omega $.}
\item{The limit $ \Omega_r := \lim_{t \rightarrow \infty} \omega_r (t,.) $ exists and is a diffeomorphism from $ {\mathcal S} $ to $ {\mathcal S} $.}
\item{For any coordinate system $ \theta_{\mathcal S} $  associated to the cover (\ref{recouvrement}), we have
$$ \theta_{\mathcal S} \circ \big( \Omega_r^{-1} \circ \omega_r \big) \big( t, {\theta}^{-1}_{\mathcal S} (\theta) \big) = \theta \ \ \mbox{mod} \ \ S^{-\tau}, \qquad \theta \in B (\theta_0, \epsilon_{\omega_0}) . $$}
\item{There exist $ \phi_r \in C^{\infty} ({\mathcal S}, \Ra) $ such that
$$ x_r (t,\omega) = t + \phi_r (\omega) \ \ \mbox{mod} \ S^{-\tau}, $$
for $ t \geq 0 $ and $ \omega \in {\mathcal S} $.}
\item{For all $ r \gg 1 $, $ N_r $ is a homeomorphism (resp. diffeomorphism) from $ [0,\infty) \times {\mathcal S} $ onto $ [r,+\infty) \times {\mathcal S} $ (resp. $(r, + \infty) \times {\mathcal S}$).}
\end{enumerate}
\end{prop}

This proposition will follow from Proposition \ref{propositiondifferee} and the following lemma on perturbations of the identity (see Appendix \ref{preuvelemmediffeo} for the proof).

\begin{lemm} \label{diffeoprocheidentite} Let $ F_{t,r} : {\mathcal S} \rightarrow {\mathcal S} $ be a family of smooth maps indexed by $ r \gg 1 $ and $ t \geq 0 $, such that, for some $ C > 0 $,
\begin{eqnarray}
 d (F_{t,r}(\omega),\omega) \leq C \scal{r}^{-\tau} , \qquad \qquad r \gg 1, \ \omega \in {\mathcal S}, \ t \geq 0  , \label{distanceabsolue}
\end{eqnarray}
and, in each chart of the cover (\ref{recouvrement}),
\begin{eqnarray}
 \big| \big| D \big(\theta_{\mathcal S} \circ F_{t,r} \circ \theta_{\mathcal S}^{-1}\big) (\theta) - I_{\Ra^{n-1}} \big| \big| \leq C \scal{r}^{-\tau} , \qquad r \gg 1, \ \theta \in B (\theta_0 , 2 \epsilon_{\omega_0}) , \ t \geq 0 \ . \label{operatornorm}
\end{eqnarray}
Then, for all $ r $ large enough and all $ t \geq 0 $, $ F_{t,r} $ is a smooth diffeomorphism on $ {\mathcal S} $.
\end{lemm}
In (\ref{operatornorm}), $ || \cdot || $ is a fixed norm on linear maps on $ \Ra^{n-1} $. Note also that $  \theta_{\mathcal S} \circ F_{t,r} \circ \theta_{\mathcal S}^{-1} $ is meaningful on $ B (\theta_0, 2 \epsilon_{\omega_0}) $, since (\ref{distanceabsolue}) implies, if $r$ is large enough, that $ F_{t,r} $ maps $ \theta_{\mathcal S}^{-1} (B (\theta_0,2 \epsilon_{\omega_0})) $ into $ \theta_{\mathcal S}^{-1} (B (\theta_0,3 \epsilon_{\omega_0})) $ which is contained in the domain of $ \theta_{\mathcal S} $ by (\ref{conditionreconnaissance}).

\bigskip

\noindent {\it Proof of Proposition \ref{propnormalflow}.} For $r$ large enough, (\ref{quantitatif}) allows to assume that $ a^{-1/2} (r,\theta) \in [1/2,3/2 ] $ hence that the initial condition $ (r,\theta,a^{-1/2},0) $ satisfies the assumption (\ref{conditionsinitiales2}). By (\ref{expressiongeodesique}), (\ref{stabilitedomaine}) and (\ref{traductioncarte}), we have then
$$ \omega \in \theta_{\mathcal S}^{-1} \big( B (\theta_0,2 \epsilon_{\omega_0}) \big) \qquad \Longrightarrow \qquad \omega_r (t, \omega) \in \theta_{\mathcal S}^{-1} \big( B (\theta_0, 3 \epsilon_{\omega_0}) \big) $$
and, by (\ref{flotpratique2}), 
$$ \big| \theta_{\nu_r}(t,\theta) - \theta \big| = \left| \frac{1}{2} \int_0^t \partial_s \theta^{s/2} (r,\theta,a^{-1/2},0) ds \right| \lesssim \scal{r}^{-\tau} , \qquad r \gg 1 , \theta \in B (\theta_0, 2 \epsilon_{\omega_0}), \ t \geq 0 . $$
This is {\it a fortiori} true if $ \theta \in B (\theta_0, \epsilon_{\omega_0}) $. So we obtain,  using (\ref{recouvrement}) and (\ref{globalisationdistance2}), that
\begin{eqnarray}
 d (\omega , \omega_r (t,\omega)) \leq C \scal{r}^{-\tau} , \qquad r \gg 1, \ \omega \in {\mathcal S}, \ t \geq 0 . \label{rappeldiffeo1}
\end{eqnarray}
Furthermore, by (\ref{flotpratique2}), we also see that $ \theta_{\nu_r} (t,.) = \theta_{\mathcal S} \circ \omega_r \big(t, \theta_{\mathcal S}^{-1}(.) \big) $ satisfies 
the condition (\ref{operatornorm}), since
\begin{eqnarray}
 \big| \partial_{\theta} \big( \theta_{\nu_r}(t,\theta) - \theta \big) \big| & = & \left| \frac{1}{2} \int_0^t \big( \partial_s \partial_{\theta} \theta^{s/2} \big) (r,\theta,a^{-1/2},0) + \big( \partial_s \partial_{\rho} \theta^{s/2} \big) (r,\theta,a^{-1/2},0) \partial_{\theta} a^{-1/2} ds \right| \nonumber \\
 & \lesssim & \scal{r}^{-\tau} , \label{rappeldiffeo2} 
\end{eqnarray}
for all $ r \gg 1 $, $ t \geq 0 $ and $ \theta \in B (\theta_0, 2 \epsilon_{\omega_0}) $\footnote{this is the interest of considering initial conditions with $ \theta \in B (\theta_0,2 \epsilon_{\omega_0}) $}. This proves the item 1.

We now consider the item 2. To prove the existence of the limit of $ \omega_r (t,.) $ as $t$ goes to infinity, it suffices to show that $ \theta_{\nu_r} (t,\theta) $ has a limit for each $ \theta \in B (\theta_0 , \epsilon_{\omega_0}) $, since we now that, by taking $r$ large enough, $ \omega_r (t,\omega) $ belongs to $ \theta_{\mathcal S}^{-1}(B (\theta_0,2 \epsilon_{\omega_0})) $ if $ \omega \in \theta_{\mathcal S}^{-1}(B (\theta_0, \epsilon_{\omega_0})) $. The existence of the limit will then follow from the integrability of $ \partial_t \theta_{\nu_r} $, which is an immediate consequence of
$$ \partial_{t} \theta_{\nu_r} (t,\theta) = \frac{1}{2} \partial_t \theta^{t/2} (r,\theta,a^{-1/2},0) = {\mathcal O} (\scal{r+t}^{-1-\tau}) $$
by (\ref{flotpratique2}). The derivatives with respect to $ \theta $ satisfy the same  bounds in time, so the limit as $ t \rightarrow \infty $ of $ \theta_{\nu_r}(t,.) $ is smooth. We can also let $t$ go to infinity in (\ref{rappeldiffeo1}) and (\ref{rappeldiffeo2}) to conclude that $ \Omega_r $ satisfies the assumptions of Lemma \ref{diffeoprocheidentite} and thus is a diffeomorphism for $r$ large enough.

To prove the item 3, we start by choosing $r$ large enough so that 
$$ \theta_{\nu_r} \left( t , B (\theta_0, \epsilon_{\omega_0 }) \right) \subset B (\theta_0, \frac{3}{2} \epsilon_{\omega_0} ) . $$
Furthermore, since $ \Omega_r $ satisfies the same bound as $ \omega_r (t,.) $ in (\ref{rappeldiffeo1}), this also holds for $ \Omega_r^{-1} $. So we may assume that
$$ \Omega_r^{-1} \left( B (\theta_0, \frac{3}{2} \epsilon_{\omega_0} ) \right) \subset B (\theta_0, 2 \epsilon_{\omega_0} ) . $$
Thus, by setting $ \Theta := \theta_{\mathcal S} \circ \Omega_r^{-1} \circ \theta_{\mathcal S}^{-1} $, it suffices to consider $ \Theta \circ \theta_{\nu_r} $. Since $ \theta = \lim_{t \rightarrow \infty} \Theta \circ \theta_{\nu_r}(t,\theta) $, we have
\begin{eqnarray*}
 \theta -  \Theta \circ \theta_{\nu_r} (t,\theta) & = & \int_t^{+ \infty} \partial_s \big( \Theta \circ \theta_{\nu_r} \big) (s,\theta) ds \\
 & = & \int_t^{+ \infty} \big( D \Theta \big) ( \theta_{\nu_r} (s,\theta) ) \cdot \partial_s \theta_{\nu_r}(s,\theta) ds  \ = \ {\mathcal O} \big( \scal{t}^{-\tau} \big)
\end{eqnarray*}
using (\ref{expressiongeodesique})  and (\ref{flotpratique2}). By differentiating this expression in $t$ and $ \theta$, we conclude that $ \Theta \circ \theta_{\nu_r} - \theta $ belongs to $ S^{-\tau} $, which is the expected result.

To prove the item 4, we observe first that the existence of $ \phi_r $ is equivalent to the existence of $ \lim_{t \rightarrow + \infty} ( x_r (t,.) - t ) $ which follows from the integrability of $ \partial_t x_r - 1 $.  This integrability in turn follows from (\ref{flotpratique1}) and (\ref{energieautomatique}) using the local expression of  $ x_r $ given by (\ref{expressiongeodesique})  and (\ref{traductioncarte}). We actually have  the following formula
\begin{eqnarray}
 x_{\nu_r} (t,\theta) & = & t - r + \int_0^t \big( \partial_s x_{\nu_r} (s,\theta) - 1 \big) ds  \nonumber \\
  & = & t + \big( \phi_r \circ \theta_{\mathcal S}^{-1} \big) (\theta) -  \int_t^{\infty} \big( \partial_s x_{\nu_r} (s,\theta) - 1 \big) ds . \label{lastterm}
\end{eqnarray}
Since $ \partial_{\theta}^{\alpha} (\partial_t x_r - t) $ is  integrable in time for any $ \alpha $, we see that $ \phi_r $ is smooth. It also follows easily from (\ref{flotpratique1}) that the last term in (\ref{lastterm}) belongs to $ S^{-\tau} $.

It remains to prove the item 5. It is convenient to denote by $ O_r (t,.) : {\mathcal S} \rightarrow {\mathcal S} $ the inverse map of $ \omega_r (t,.) $. Note that since $ \omega_r $ is smooth on $ [0,\infty ) \times {\mathcal S} $,  so is the map $ O_r : (t,\omega) \mapsto O_r (t,\omega)$. Therefore, the map
$$ M_r : \ (t,\omega) \mapsto (t, \omega_r (t,\omega))  $$ 
is a homeomorphism from $ [0,\infty ) \times {\mathcal S} $  onto itself with inverse $ (t,\omega) \mapsto (t, O_r (t,\omega)) $. It is also obviously a diffeomorphism on the interior. It is thus sufficient to prove the result for the map $ P_r := N_r \circ M_r^{-1} $ instead of $N_r $. Notice that $ P_r $ has the following simpler form
$$ P_r (t,\omega) = \big(  x_r (t,O_r (t,\omega)) , \omega \big) . $$
This map is smooth up to $t = 0 $ and it is thus not hard to see that the conclusion would be a consequence of the fact that, for each $ \omega \in {\mathcal S} $, the map
$$ t \mapsto X_{r,\omega}(t) :=  x_r (t,O_r (t,\omega)) $$
is a bijection from $ [0,\infty ) $ onto $ [r,\infty ) $. Clearly, if $t=0$ we have $ X_{r,\omega} (0) = r $, so it is sufficient to show that
\begin{eqnarray}
 |\partial_t X_{r,\omega} (t) - 1 | \leq 1/2 , \label{deriveeenrjuste}
\end{eqnarray}
for $r$ large enough and $ t \geq 0 $. Using (\ref{flotpratique2}) and (\ref{rappeldiffeo2}), it is not hard to see that $ \partial_t \theta_{\mathcal S} \circ O_r (t,.) $ is of order $ \scal{r}^{-\tau} $ which, together with (\ref{flotpratique1}), implies (\ref{deriveeenrjuste}) and completes the proof. \finpreuve

\bigskip

\subsection{Proof of Theorem \ref{theoremecoordonnees}} \label{subsection2point2}
Item 1 follows from Proposition \ref{propositiondifferee} and (\ref{expressiongeodesique}). The item 2 is the item 5 of Proposition \ref{propnormalflow}. We now prove the item 3. If $ \theta_{\mathcal S} = ( \theta_1 , \ldots , \theta_{n-1} ) $ are coordinates on $ {\mathcal S} $, then $ (t, \theta_1 , \ldots , \theta_{n-1} ) $ are coordinates on $ (0,\infty) \times {\mathcal S} $ and
$$ \overline{t} := t \circ N_r^{-1}, \qquad \overline{\theta}_j := \theta_j \circ N_r^{-1}, \qquad j = 1, \ldots , n-1 , $$
are coordinates on $ {\mathcal M} $ which we work with. It is useful to note, by standard properties of the local normal flow, that $ N_r $ is smooth up to $ t= 0 $ and $ N_r^{-1} $ up to $ x =r $. In particular, this allows us to use the fact that the vector fields $ \partial / \partial t $, $ \partial / \partial \theta_j $, $ \partial / \partial \overline{t} $ and $ \partial / \partial \overline{\theta}_j $ are defined up to the boundary. We  show first that
\begin{eqnarray}
 N_r^* {\mathbf G}_{(t,\omega)} \left( \frac{\partial}{\partial t}, \frac{\partial}{\partial t} \right) = 1 , \qquad N_r^* {\mathbf G}_{(t,\omega)} \left( \frac{\partial}{\partial t}, \frac{\partial}{\partial \theta_j} \right) = 0 . \label{warpedfinal}  
\end{eqnarray}
To that end, we observe on one hand that
\begin{eqnarray}
\frac{\partial}{\partial \overline{t}}{|_{N_r(t,\omega)}} = d N_r  \left( \frac{\partial}{\partial t} {|_{(t,\omega)}} \right), \qquad \frac{\partial}{\partial \overline{\theta}_j}{|_{N_r(t,\omega)}} = d N_r  \left( \frac{\partial}{\partial \theta_j} {|_{(t,\omega)}} \right) ,
\end{eqnarray}
and, on the other hand that
\begin{eqnarray}
 \frac{\partial}{\partial \overline{t}}{|_{N_r(t,\omega)}} =  \frac{d}{dt} N_r (t,\omega) , \label{geodesiqueintrinseque}
\end{eqnarray}
which is the tangent vector to the geodesic $ \exp_{(r,\omega)} (t \nu_r) $. In particular, at $ t= 0$, the vector field in (\ref{geodesiqueintrinseque}) is $ \nu_r $ so  (\ref{warpedfinal}) is true for $t = 0 $. It then suffices to show that the left hand sides in (\ref{warpedfinal}) are constant with respect to $ t $.
Using the standard properties of the Levi-Civita connection $ \nabla $ and (\ref{geodesiqueintrinseque})
\begin{eqnarray*}
\frac{\partial}{\partial t} N_r^* {\mathbf G}_{(t,\omega)} \left( \frac{\partial}{\partial t}, \frac{\partial}{\partial t} \right) & =& \frac{\partial}{\partial t} \left(
{\mathbf G} \left( \frac{\partial}{\partial \overline{t}}, \frac{\partial}{\partial \overline{t}} \right) \right)_{| N_r (t,\omega)} \\
& = & \frac{\partial}{\partial \overline{t} }  
{\mathbf G} \left( \frac{\partial}{\partial \overline{t}}, \frac{\partial}{\partial \overline{t}} \right) \\
& = & 2 {\mathbf G} \left( \nabla_{\frac{\partial}{\partial \overline{t}}} \frac{\partial}{\partial \overline{t}} , \frac{\partial}{\partial \overline{t}} \right) \ = \ 0 ,
\end{eqnarray*}
where, in the last two lines, we dropped the evalutation at $ N_r (t,\omega) $ from the notation for simplicity.
This yields the first equality of (\ref{warpedfinal}) for all $ t \geq 0$. 
For the second equality, we compute similarly
\begin{eqnarray*}
\frac{\partial}{\partial t} N_r^* {\mathbf G}_{(t,\omega)} \left( \frac{\partial}{\partial t}, \frac{\partial}{\partial \theta_j} \right) & =& \frac{\partial}{\partial t} \left(
{\mathbf G} \left( \frac{\partial}{\partial \overline{t}}, \frac{\partial}{\partial \overline{\theta}_j} \right) \right)_{| N_r (t,\omega)} \\
& = & \frac{\partial}{\partial \overline{t} }  
{\mathbf G} \left( \frac{\partial}{\partial \overline{t}}, \frac{\partial}{\partial \overline{\theta}_j} \right) \\
& = &  {\mathbf G} \left( \nabla_{\frac{\partial}{\partial \overline{t}}} \frac{\partial}{\partial \overline{t}} , \frac{\partial}{\partial \overline{\theta}_j} \right) + 
{\mathbf G} \left(  \frac{\partial}{\partial \overline{t}}, \nabla_{\frac{\partial}{\partial \overline{t}}} \frac{\partial}{\partial \overline{\theta}_j}  \right) .
\end{eqnarray*}
Here, using that the Levi-Civita connection is torsion free, we have
\begin{eqnarray*}
 {\mathbf G} \left(  \frac{\partial}{\partial \overline{t}}, \nabla_{\frac{\partial}{\partial \overline{t}}} \frac{\partial}{\partial \overline{\theta}_j}  \right) & = &
{\mathbf G} \left(  \frac{\partial}{\partial \overline{t}}, \nabla_{\frac{\partial}{\partial \overline{\theta}_j}} \frac{\partial}{\partial \overline{t}} + \left[ \frac{\partial}{\partial \overline{t}} , \frac{\partial}{ \partial \overline{\theta}_j } \right] \right) \\
& = & \frac{1}{2} \frac{\partial}{\partial \overline{\theta}_j} {\mathbf G} \left( \frac{\partial}{\partial \overline{t}}, \frac{\partial}{\partial \overline{t}} \right) \ = \ 0 ,
\end{eqnarray*}
since the Lie bracket in the first line vanishes and since, in the second line,  we are differentiating a constant function. 
This completes the proof of (\ref{warpedfinal}).

To determine $  N_r^* {\mathbf G} \left( \frac{\partial }{\partial \theta_i} ,  \frac{\partial }{\partial \theta_j} \right) $ it suffices to compute the last $ n-1 $ columns and rows of the following block matrix decomposition of $ N_r^* {\mathbf G} $ in local coordinates,  
$$ \left( \begin{matrix}
\partial x_{\nu_r} / \partial t & \partial x_{\nu_r} / \partial \theta     \\
\partial \theta_{\nu_r} / \partial t &  \partial \theta_{\nu_r} / \partial \theta 
\end{matrix}
 \right)^T \left( \begin{matrix}
1 & 0     \\
0 & w    
\end{matrix}
 \right)^{-1} \left( \begin{matrix}
{\mathbf a} & {\mathbf b}^T    \\
{\mathbf b} & {\mathbf g}   
\end{matrix}
 \right) \left( \begin{matrix}
1 & 0     \\
0 & w  \end{matrix}
 \right)^{-1} \left( \begin{matrix}
\partial x_{\nu_r} / \partial t & \partial x_{\nu_r} / \partial \theta     \\
\partial \theta_{\nu_r} / \partial t &  \partial \theta_{\nu_r} / \partial \theta 
\end{matrix}
 \right) ,  $$
 where $ {\mathbf a} $, $ {\mathbf b} $, $ {\mathbf g} $ and $ w $ are evaluated at $ (x_{\nu_r},\theta_{\nu_r})(t,\theta) $. After a simple calculation, the matrix is
\begin{eqnarray}
 w^{-2} \left\{ w^2 {\mathbf a} \frac{\partial x_{\nu_r}}{\partial \theta}^T \frac{\partial x_{\nu_r}}{\partial \theta} + w \left(  \frac{\partial x_{\nu_r}}{\partial \theta}^T {\mathbf b}^T \frac{\partial \theta_{\nu_r}}{\partial \theta} + \frac{\partial \theta_{\nu_r}}{\partial \theta}^T {\mathbf b} \frac{\partial x_{\nu_r}}{\partial \theta} \right) +  \frac{\partial \theta_{\nu_r}}{\partial \theta}^T {\mathbf g} \frac{\partial \theta_{\nu_r}}{\partial \theta} \right\} . \label{nouvellemetriquealinfini}
\end{eqnarray}
By (\ref{quantitatif}), (\ref{conditionexposant}) and Proposition \ref{propnormalflow}, the matrix (of the metric) inside $ \{ \cdots \} $ is of the form
\begin{eqnarray}
 S^{- 1 - \tau} + S^{- 1 - \tau} +  \big( \theta_{\mathcal S}^{-1} \big)^* \Omega_{r}^* \overline{\mathbf g} + S^{-\tau} ,
\label{facteurconformeconcret} 
\end{eqnarray}
where, for the last two terms, we used that $ \theta_{\nu_r}(t,.) = \theta_{\mathcal S} \circ \Omega_r \circ \theta_{\mathcal S}^{-1} + S^{-\tau} $  as well as the fact that $ \theta_{\mathcal S}^* {\mathbf g} = \overline{\mathbf g} + S^{-\tau} $.
On the other hand, using the second condition of (\ref{hypothesew}) and (\ref{limitew}), we have $ w^{\prime} / w - \kappa  \in S^{- \varepsilon} $, from which it follows that
\begin{eqnarray}
 w (t + b) = w (t)  e^{ \kappa b } \exp \left( \int_t^{t+b} \sigma_{-\varepsilon}(s)ds \right) ,
\label{influenceconforme}
\end{eqnarray}
for some $ \sigma_{-\varepsilon} \in S^{- \varepsilon} $. This identity and the item 4 of Proposition \ref{propnormalflow} imply that
\begin{eqnarray}
 w (x_{\nu_r}(t,\theta)) = w (t) e^{ \kappa \left( \phi_r \circ \theta_{\mathcal S}^{-1} \right)(\theta)} \left( 1 + S^{- \min (\varepsilon,\tau)} \right) . 
\nonumber
\end{eqnarray}
Combining this identity and (\ref{facteurconformeconcret}) completes the proof of the item 3 of Theorem \ref{theoremecoordonnees}. 
\finpreuve

\bigskip

\noindent {\bf Justification of example 1.} Using the item 4 of Proposition \ref{propnormalflow}, we see that the term $ w (x_r)^{-2} $ in front of (\ref{nouvellemetriquealinfini}) is of the form
$$ w(x_r)^{-2} = (t + \phi_r + S^{-\tau})^2 = t^{2} \left(1 + 2 \phi_r t^{-1} + o(t^{-1})\right) ,$$
which proves (\ref{expansionexample1}).
\bigskip

\noindent {\bf Justification of example 3.} In this case, (\ref{influenceconforme}) reads explicitly
\begin{eqnarray}
 w (t+b) & = & w (t) e^{-b} \exp \left( - \int_0^b \beta (t+u)^{\beta - 1} du \right)  \nonumber \\
  & = & w (t) e^{-b} \big( 1 - \beta b t^{\beta-1} + o (t^{\beta-1}) \big) , \label{thisexpansion}
\end{eqnarray}
where $ o (t^{\beta-1}) $ is uniform with respect to $ b$ as long as $b$ remains in a compact set. Using again the item 4 of Proposition \ref{propnormalflow} to write $ x_r $ as $ t +  b$, (\ref{thisexpansion}) combined with (\ref{nouvellemetriquealinfini}) and (\ref{facteurconformeconcret}) implies (\ref{expansionexample3}).

\section{Proof of Proposition \ref{propositiondifferee}} \label{subsectionproof}
\setcounter{equation}{0}
The proof will be reduced to the analysis of hamiltonians globally defined on $ \Ra^{2n} $. Indeed, by possibly increasing $ R $ and by (\ref{hypothesew}), we may assume that $ w  $ is defined on $ \Ra $ and belongs to $ S^{-\lambda} (\Ra) $. Also, by (\ref{quantitatif}), we can modify the coefficients of $p$ on $ B (\theta_0,4 \epsilon_{\omega_0}) \setminus B (\theta_0, 3 \epsilon_{\omega_0}) $ so that
\begin{eqnarray}
 a -1 \in S^{- \min (\mu,2 \nu)} (\Ra \times \Ra^{n-1}), \qquad b \in S^{-\nu} (\Ra \times \Ra^{n-1}) , \qquad g - \bar{\mathbf g}^{-1} \in S^{-\tau} (\Ra \times \Ra^{n-1}) , 
 \label{hamRn}
\end{eqnarray}
for some positive definite matrix $ \bar{\mathbf g}^{-1}  $ defined on $ \Ra^{n-1} $ with $ C^{\infty}_b $ coefficients, such that $ \bar{\mathbf g}^{-1} (\theta) \geq C > 0 $ for all $\theta$ and which coincides with the original $ \bar{\mathbf g}^{-1} $ on $ B (\theta_0, 3 \epsilon_{ \omega_0}) $. Then, we keep the notation $p$ for the symbol
\begin{eqnarray}
 p (x,\theta,\rho,\eta) = a (x,\theta)\rho^2 + 2 w(x) \rho b (x,\theta) \cdot \eta + w(x)^2 \eta \cdot g (x,\theta) \eta  , 
\label{formesymboltotalbis}
\end{eqnarray}
which coincides with the principal symbol of the Laplacian on $ (R,+\infty) \times B (\theta_0,3 \epsilon_{\omega_0}) \times \Ra^n $. We may  assume that
for some $ C_0 \geq 1 $,
\begin{eqnarray}
C_0^{-1} (\rho^2 + w(x)^2 |\eta|^2) \leq  p(x,\theta,\rho,\eta) \leq C_0 \big( \rho^2 + w(x)^2 |\eta|^2 \big) , \label{ellipticiteuniforme}
\end{eqnarray}
everywhere on $ \Ra^{2n} $. 

We consider $ \big( x^t,\theta^t, \rho^t , \eta^t \big) $, the hamiltonian flow  of $p$
with initial condition $ \big(x,\theta,\rho ,\eta \big) $ at $ t= 0 $.

\begin{prop} \label{premierepartie}  Assume (\ref{hypothesew}), (\ref{quantitatif}) and (\ref{conditionexposant}).
Then, for all $ M > 1 $, there exists $  X_1 >0 $ such that, for all
\begin{eqnarray}
 x \geq X_1, \qquad \theta \in \Ra^{n-1}, \qquad \rho \in \big[ M^{-1}, M \big], \qquad |\eta| \leq M ,  \label{conditionsinitiales}
\end{eqnarray}
the hamiltonian flow of $p$ is defined for all $  t \geq 0 $ and satisfies
\begin{eqnarray}
\left\{ \begin{array}{r c l} \big| x^t - x - 2 t \rho^t \big| & \lesssim & \scal{x}^{-\tau} , 
\\
\big| \theta^t - \theta \big| & \lesssim & \scal{x}^{-\tau}, 
\\
\big| \rho^t \big| & \lesssim & 1 ,  
\\
\big| \eta^t \big| & \lesssim & 1 , 
\end{array} \right. \label{bornexyz}
\end{eqnarray}
where $ p = p (x,\theta,\rho,\eta) $. Furthermore, for all $t \geq 0$
\begin{eqnarray}
x^t & \geq & x + \frac{t}{M},  \label{stabilitedomaine2} \\
\rho^t & \gtrsim  & 1 , \label{minorationrho} \\
|\rho^t - p^{1/2}| & \lesssim & \scal{x + t}^{-1-\tau} . \label{convergencequantitative}
\end{eqnarray}
\end{prop}

Notice that (\ref{convergencequantitative}) implies that
\begin{eqnarray}
\lim_{t \rightarrow + \infty} \rho^t = p^{1/2} , \label{momentasymptotique}
\end{eqnarray}
and also that, in the left hand side of the first estimate of (\ref{bornexyz}), $ 2 t \rho^t $ could be replaced by $ 2 t p^{1/2} $.

\bigskip

\noindent {\it Proof.} By boundedness of $ w $ and $ w^{\prime} $, we have
\begin{eqnarray}
 p \big( x,\theta,\rho,\eta \big) \leq C^{\prime}_0 , \qquad \mbox{for} \ \ |\rho| \leq M, \ |\eta| \leq M , \label{pourenergieconservee}
\end{eqnarray}
with $ C^{\prime}_0 $ depending only on $ C_0 $ and $ M $. On the other hand, by (\ref{conditionexposant}) and (\ref{hamRn}), we have  
\begin{eqnarray}
\left\{ \begin{array}{r c l}
\big| \partial_{\rho} p - 2 \rho \big| & \leq & C_1 \scal{x}^{-1-\tau} \big( |\rho| + |\eta| \big) , 
\\
\big| \partial_{\eta} p \big| & \leq & C_2 \scal{x}^{-1-\tau} \big( |\rho| + |\eta| \big) , 
\\
\big| \partial_x p \big| & \leq & C_3  \scal{x}^{-2-\tau} \big( \rho^2 + |\eta|^2 \big) , 
\\
\big| \partial_{\theta} p \big| & \leq & C_4 \scal{x}^{-1-\tau} \big( \rho^2 + |\eta|^2 \big) , 
\end{array}
\right. \qquad \mbox{on} \ \Ra^{2n} , \label{unxyz}
\end{eqnarray}
using that 
$$ \min (\mu , 2 \nu) \geq 1 + \tau, \qquad \lambda + \nu \geq 1 + \tau , \qquad 2 \lambda \geq 1 + \tau .  $$
 Given $ (x,\theta,\rho,\eta) $ satisfying (\ref{conditionsinitiales}),
denote by $ [0,T_+) $ the domain of the maximal solution. We shall prove that $ T_+ = + \infty $ and that
\begin{eqnarray}
x^t \geq x + \frac{t}{M}, \qquad |\eta^t| \leq 2 M , \label{inegaliteI}
\end{eqnarray}
for all $t \in [0,T_+ )$. Introduce the set
$$ I := \left\{ T \in [0,T_+ ) \ | \ (\ref{inegaliteI}) \ \ \mbox{holds  on} \  [0,T] \right\} .  $$
This is obviously an interval containing $ 0 $ and we set $ T_{++} = \sup I $, which is clearly positive. 
 Using (\ref{pourenergieconservee}), the conservation of energy and (\ref{ellipticiteuniforme}), we obtain a bound 
 $$ |\rho^t | \leq \big( C_0 C_0^{\prime} \big)^{1/2} $$ along the flow and see that there exist $ C^{\prime}_1 , C^{\prime}_3 , C^{\prime}_4 $ depending only on $ C_1,C_3,C_4 $ and $ M $ such that 
\begin{eqnarray}
\big| \dot{x}^s - 2 \rho^s \big| & \leq & C_1^{\prime} \scal{x^s}^{-1-\tau}  , \nonumber  \\
\big| \dot{\rho}^s \big| & \leq & C_3^{\prime}  \scal{x^s}^{-2-\tau} , \nonumber \\
\big| \dot{\eta}^s \big| & \leq & C_4^{\prime} \scal{x^s}^{-1-\tau}  , \nonumber
\end{eqnarray}
for all $ s \in I $. Thus, if one chooses $ X_1 $ large enough so that
\begin{eqnarray}
 C^{\prime}_3  \int_0^{\infty} \left\langle X_1 + \frac{s}{M} \right\rangle^{-1-\tau}d s & < & \frac{1}{4M}, \nonumber
 \\
 C^{\prime}_4  \int_0^{\infty} \left\langle X_1 + \frac{s}{M} \right\rangle^{-1-\tau}d s & < & \frac{M}{4}, \nonumber \\
  C^{\prime}_1 \scal{X_1}^{-\tau} < \frac{1}{4M} , \nonumber
\end{eqnarray}
then, for all $t \in I$,
\begin{eqnarray*}
 \dot{x}^t   \geq  2 \rho^t - \frac{1}{4M} , \qquad
\big| {\rho}^t - \rho \big|  \leq  \frac{1}{4M} , \qquad
\big| \eta^t - \eta \big|  \leq  \frac{M}{4}  .
\end{eqnarray*}
Using (\ref{conditionsinitiales}), this implies clearly that, for all $t \in I$,
\begin{eqnarray*}
\big| \eta^t  \big|  \leq  \frac{5M}{4}  , \qquad
\rho^t  \geq  \frac{3}{4M} , \qquad
x^t  \geq  x + \frac{5}{4} \frac{t}{M} ,
\end{eqnarray*}
yielding a contradiction with the fact that $ T_{++} < T_+ $ (one could otherwise obtain (\ref{inegaliteI}) beyond $ T_{++} $). Thus $ T_{++} = T_+ $ and $ T_+ = + \infty $, since   (\ref{unxyz}) and (\ref{inegaliteI})  imply that the flow cannot blow up in finite time. We have thus shown the completness of flow on $ [ 0 , + \infty ) $ as well as the third and fourth estimates of (\ref{bornexyz}), (\ref{stabilitedomaine2}) and (\ref{minorationrho}). In particular, using that $ x^t \rightarrow \infty $ as $ t \rightarrow + \infty $, we also deduce (\ref{momentasymptotique})  from the conservation of energy and  the positivity of $ \rho^t $. Integrating $ \dot{\rho}^s $ for $s \in [t,\infty) $, we obtain the quantitative bound (\ref{convergencequantitative}), using the third estimate of (\ref{unxyz}) and (\ref{inegaliteI}).
  It remains to prove the first two estimates of (\ref{bornexyz}). For the first one, it suffices to observe that 
 $$  \big| \partial_t (x^t - x -2 t \rho^t) \big| = \big| \dot{x}^t - 2 \rho^t - 2 t \dot{\rho}^t \big| \lesssim \scal{x + t}^{-1-\tau} , $$
 using the third estimate of (\ref{bornexyz}),  the first and third estimates of (\ref{unxyz}) and (\ref{inegaliteI}).   The second one is obtained similarly from the second estimate of (\ref{unxyz}). \finpreuve

\bigskip

\noindent {\bf Remark.} As one can see from this proof, the completness of the flow as well as the estimates  (\ref{bornexyz}) (3rd and 4th) to (\ref{minorationrho}) could be obtained even if we only had $ -\tau $ and $ - 1 - \tau $ rather than $ -1- \tau $ and $ -2- \tau $ in the first and third lines of (\ref{unxyz}) respectively. Furthermore, in this case we also would have a lower bound similar to (\ref{stabilitedomaine2}).
The powers $ - 1 - \tau $ and $ - 2 -\tau $ play a role only when we prove the first estimate of (\ref{bornexyz}).

\bigskip

For future reference, we note here the following elementary fact. Assuming that the initial conditions satisfy (\ref{conditionsinitiales}) with $ X_1$ large enough, we can  freely modify the Hamiltonian vector field of $p$ for $ |\rho| + |\eta|  $ large ({\it e.g.} cutoff) by conservation of energy. More precisely, using the last two estimates of (\ref{bornexyz}), we work on a domain where we can assume that the Hamilton equations (\ref{hamiltonequation}) read
\begin{eqnarray}
 \begin{cases}
\dot{x}^t & = \ 2 \rho^t + a_1 \big(x^t,\theta^t,\rho^t,\eta^t  \big) \ = \ a_0 \big(x^t,\theta^t,\rho^t,\eta^t  \big), \\
\dot{\theta}^t & = \ a_2 \big(x^t,\theta^t,\rho^t,\eta^t  \big) , \\
\dot{\rho}^t & = \ a_3 \big(x^t,\theta^t,\rho^t,\eta^t  \big),  \\
\dot{\eta}^t & = \  a_4 \big(x^t,\theta^t,\rho^t,\eta^t  \big),
\end{cases}
   \label{Hamilton}
\end{eqnarray}
with
$$ a_1,a_2,a_4 \in S^{-\tau-1} , \qquad a_3 \in S^{-\tau-2}, \qquad a_0 \in S^0 . $$
This remark will be useful below. In the next proposition, we recall that $ \partial^{\gamma} = \partial_{x}^k \partial_{\theta}^{\alpha} \partial_{\rho}^l \partial_{\eta}^{\beta} $.
\begin{prop} \label{propositionb} Assume (\ref{hypothesew}), (\ref{quantitatif}) and (\ref{conditionexposant}). Then, for all $ M > 0 $, there exists $ X_1 > 0 $ such that, on the domain defined by (\ref{conditionsinitiales}), we have
\begin{eqnarray}
\left\{ \begin{array}{ r c l}
\big| \partial^{\gamma} (x^t-x-2t \rho^t) \big| & \lesssim & \scal{x}^{-\tau} , 
\\
\big| \partial^{\gamma} (\theta^t-\theta) \big| & \lesssim & \scal{x}^{-\tau} , 
\\
\big| \partial^{\gamma} (\rho^t-\rho \big) \big| & \lesssim & \scal{x}^{-\tau-1} , \\
\big| \partial^{\gamma} (\eta^t-\eta) \big| & \lesssim & \scal{x}^{-\tau} , 
\end{array} \right. \label{sansj}
\end{eqnarray}
and, for $ j \geq 1 $,
\begin{eqnarray}
\left\{ \begin{array}{ r c l}
\big| \partial_t^j \partial^{\gamma} (x^t-x-2t \rho^t) \big| & \lesssim & \scal{x+t}^{-\tau-j} , 
\\
\big| \partial_t^j \partial^{\gamma} (\theta^t-\theta) \big| & \lesssim & \scal{x+t}^{-\tau-j} ,
 \\
\big| \partial_t^j \partial^{\gamma} (\rho^t-\rho \big) \big| & \lesssim & \scal{x+t}^{-\tau-1-j} , 
\\
\big| \partial_t^j \partial^{\gamma} (\eta^t-\eta) \big| & \lesssim & \scal{x+t}^{-\tau-j} . 
\end{array} \right. \label{avecj}
\end{eqnarray}
\end{prop}

Notice that $ \rho $ may be omitted in the third line of (\ref{avecj}) or even be replaced by $ p^ {1/2} $. From this remark, we obtain the additional useful estimates, for $ j \geq 0$,
\begin{eqnarray}
\big| \partial_t^j \partial^{\gamma} (\rho^t-p^{1/2} \big) \big| & \lesssim & \scal{x+t}^{-\tau-1-j} . \label{aussiimportant}
\end{eqnarray}



\bigskip

\noindent {\it Proof.}   Let us introduce
$$ u^t := x^t  - 2 t \rho^t, \qquad  \Phi^t = \big( u^t,\theta^t,\rho^t,\eta^t \big) . $$
Clearly, (\ref{sansj}) follows by integration in time of (\ref{avecj})  since $ u^t -x $, $ \theta^t- \theta $, $ \rho^t - \rho $ and $ \eta^t - \eta $ vanish at $t=0$. It is thus sufficient to prove (\ref{avecj}), which we consider now. Using the identity
$$ \dot{u}^t = \dot{x}^t - 2 \rho^t - 2 t \dot{\rho}^t  ,$$ 
and  (\ref{Hamilton}), one checks that $ \Phi^t $ satisfies an ODE of the form
\begin{eqnarray}
 \begin{cases}
\dot{u}^t & = \ (b_1 + t \widetilde{b}_1) \big(x^t,\theta^t,\rho^t,\eta^t  \big),  \\
\dot{y}^t & = \ b_2 \big(x^t,\theta^t,\rho^t,\eta^t  \big),  \\
\dot{\rho}^t & = \ \widetilde{b}_3 \big(x^t,\theta^t,\rho^t,\eta^t  \big),  \\
\dot{\eta}^t & = \  b_4 \big(x^t,\theta^t,\rho^t,\eta^t  \big),
\end{cases}
   \label{systemereduit}
\end{eqnarray}
with
$$ b_1 , b_2 , b_4 \in S^{-\tau-1}, \qquad  \widetilde{b}_1 , \widetilde{b}_3 \in S^{-\tau-2} . $$
Independently,  (\ref{Hamilton}) again and a simple induction on $j$ show that
\begin{eqnarray}
 a \in S^m \qquad \Longrightarrow \qquad \partial_t^j a \big( x^t,\theta^t,\rho^t, \eta^t \big) = \widetilde{a} \big( x^t,\theta^t,\rho^t, \eta^t \big) \qquad \mbox{for some} \ \widetilde{a} \in S^{m-j} . \label{implicationpropre}
\end{eqnarray}
Assume for a while that we have proved the bounds
\begin{eqnarray}
 |\partial^{\gamma} \Phi^t| \leq C_{\gamma} , \qquad |\gamma| \geq 1 . \label{apriorifirst}
\end{eqnarray}
Then, for $ |\gamma| \geq 1 $, 
\begin{eqnarray}
 |\partial^{\gamma}x^t| \lesssim \scal{t}, \qquad |\partial^{\gamma}\theta^t| + |\partial^{\gamma} \rho^t| + |\partial^{\gamma} \eta^t| \lesssim 1 ,
 \label{apriorisecond}
\end{eqnarray}
and let us show how it leads to the result. By applying $ \partial_t^{j-1} $ to (\ref{systemereduit}) and using (\ref{implicationpropre}), we see first that
\begin{eqnarray}
 \begin{cases}
\partial_t^j u^t & = \ (c_1 + t \widetilde{c}_1) \big(x^t,\theta^t,\rho^t,\eta^t  \big),  \\
\partial_t^j \theta^t & = \ c_2 \big(x^t,\theta^t,\rho^t,\eta^t  \big),  \\
\partial_t^j \rho^t & = \ \widetilde{c}_3 \big(x^t,\theta^t,\rho^t,\eta^t  \big),  \\
\partial_t^j \eta^t & = \  c_4 \big(x^t,\theta^t,\rho^t,\eta^t  \big),
\end{cases} \label{reduitderive}
\end{eqnarray}
with
$$ c_1 , c_2 , c_4 \in S^{-\tau-j}, \qquad  \widetilde{c}_1 , \widetilde{c}_3 \in S^{-\tau-1-j} . $$
On the other hand, the Fa\`a Di Bruno formula (see for instance \cite{FaaDiBruno}) yields
\begin{eqnarray}
 \partial^{\gamma} \left( a (x^t,\theta^t,\rho^t,\eta^t) \right) \ = \ \partial_x a  \partial^{\gamma}x^t + \partial_{\theta} a \partial^\gamma \theta^t + \partial_{\rho}a  \partial^{\gamma} \rho^t +  \partial_{\eta}a   \partial^{\gamma} \eta^t  + \nonumber \\ \mbox{linear combination of} \ \ 
  (\partial_x^k \partial_{\theta}^{\alpha} \partial_{\rho}^l \partial_{\eta}^{\beta}a ) \prod_{1 \leq i \leq k }  \partial^{\gamma^x_i} x^t \prod_{\delta, i}  \partial^{\gamma^{\theta_{\delta}}_i} \theta_{\delta}^t \prod_{i}  \partial^{\gamma^{\rho}_i}  \rho^t \prod_{\delta,i} \ \partial^{\gamma^{\eta_{\delta}}_i}  \eta^t_{\delta} , \label{FaaDiBruno}
\end{eqnarray}
where all derivatives in the products of the second line are of striclty smaller order than $ |\gamma| $ and satisfy
$$ \sum_i \gamma^x_i  + \sum_{\delta,i} \gamma^{\theta_{\delta}}_i  + \sum_i \gamma^{\rho}_i  + \sum_{\delta,i} \gamma^{\eta_{\delta}}_i  = \gamma , $$
and where all derivatives of $a$ are of course evaluated at $ (x^t,\theta^t,\rho^t,\eta^t) $. If $a \in S^m $, using (\ref{stabilitedomaine2}),
 we deduce from (\ref{apriorisecond}) and (\ref{FaaDiBruno}) that
\begin{eqnarray*}
  \big| \partial^{\gamma} \left( a (x^t,\theta^t,\rho^t,\eta^t) \right) \big| & \lesssim & \scal{x+t}^{m-1} \scal{t} + \scal{x+t}^{m} + \sum_{k \leq |\gamma|}
\scal{x+t}^{m-k} \scal{t}^k  , \\
& \lesssim & \scal{x+t}^{m} .
 \end{eqnarray*}
Therefore, by applying $ \partial^{\gamma} $ to (\ref{reduitderive}), (\ref{avecj}) is a straightforward consequence of (\ref{apriorifirst}). It thus remains to prove (\ref{apriorifirst}), which we do now by induction on $ |\gamma| $. By (\ref{reduitderive}), we can introduce
$$ B_t = B + t \widetilde{B} , \qquad B \in S^{-\tau}, \ \widetilde{B} \in S^{-\tau-1} , $$
which are $ \Ra^{2n} $ valued so that
\begin{eqnarray}
 \dot{\Phi}^t = B_t \big( x^t,\theta^t,\rho^t,\eta^t \big). \label{equationsimple}
\end{eqnarray}
By  applying $ \partial^{\gamma} $ to this equation  (with $ |\gamma|=1 $ first)  and using that
$$ \big|\partial^{\gamma} x^t \big| \lesssim \scal{t} \big| \partial^{\gamma} X^t \big|, \qquad \big|\partial^{\gamma} \theta^t \big| + \big|\partial^{\gamma} \rho^t \big| + \big|\partial^{\gamma} \eta^t \big| \lesssim \big| \partial^{\gamma} X^t \big| , $$
we obtain
$$ |\partial^{\gamma} \Phi^t| \lesssim |\partial^{\gamma} \Phi^0| + \int_0^t \scal{x+s}^{-2-\tau}\scal{s}|\partial^{\gamma} \Phi^s| +  \scal{x+s}^{-1-\tau} |\partial^{\gamma} \Phi^s|  ds $$
using also (\ref{stabilitedomaine2}). By the Gronwall Lemma, this yields (\ref{apriorifirst}) for $ |\gamma|=1 $. Then, assuming $ |\gamma| \geq 2 $ and that (\ref{apriorifirst})
 has been proved  for lower orders, we obtain 
$$ |\partial^{\gamma} \Phi^t| \lesssim  \int_0^t \scal{x+s}^{-2-\tau}\scal{s}|\partial^{\gamma}\Phi^s| +  \scal{x+s}^{-1-\tau} |\partial^{\gamma}\Phi^s|  ds 
+ \sum_{k = 0}^{|\gamma|} \int_0^t \scal{x+s}^{-1-\tau-k} \scal{s}^k ds , $$
by applying $ \partial^{\gamma} $ to the equation (\ref{equationsimple}) and using (\ref{FaaDiBruno}). Then (\ref{apriorifirst}) follows from the Gronwall Lemma. The proof is complete. \finpreuve

\bigskip

\noindent {\bf Proof of Proposition \ref{propositiondifferee}.} The localization properties in (\ref{stabilitedomaine}) follow from (the second line of) (\ref{bornexyz}) and (\ref{stabilitedomaine2}). Note in particular that within the domain $ (X_1,\infty) \times B (\theta_0, 3 \epsilon_{\omega_0}) \times \Ra^n $ (with $ X_1 \gg 1 $), the hamiltonian flow of the globally defined hamiltonian $p$ in (\ref{formesymboltotalbis}) does indeed represent the geodesic flow in a chart. The estimates (\ref{flotpratique1}) and (\ref{flotpratique2}) follow directly from (\ref{avecj}). \finpreuve

\appendix

\section{Proof of Lemma \ref{diffeoprocheidentite}} \label{preuvelemmediffeo}
\setcounter{equation}{0}
Let us prove first that $ F_{t,r} $ is injective for $ r $ large enough. Assume that $ \omega , \omega^{\prime} \in {\mathcal S} $ satisfy $ F_{t,r} (\omega) = F_{t,r} (\omega^{\prime}) $. Then, by the triangle inequality
$$ d (\omega,\omega^{\prime}) \leq d (\omega , F_{t,r}(\omega)) + d (F_{t,r}(\omega),F_{t,r}(\omega^{\prime})) + d (F_{t,r}(\omega^{\prime}),\omega^{\prime}) \leq 2 C \scal{r}^{-\tau} . $$
For $r$ large enough, we can thus insure that if $ \omega \in \theta_{\mathcal S}^{-1} (B (\theta_0, \epsilon_{\omega_0})) $ then $ \omega^{\prime}
 \in \theta_{\mathcal S}^{-1} (B (\theta_0, 2 \epsilon_{\omega_0})) $. In particular, they belong to the same coordinate patch so we can  consider $ \theta := \theta_{\mathcal S} (\omega) $ and $ \theta^{\prime} := \theta_{\mathcal S} (\omega^{\prime}) $. Furthermore, using that $  \theta_{\mathcal S} \circ F_{t,r} (\omega) =
 \theta_{\mathcal S} \circ F_{t,r}(\omega^{\prime})  $, we have
\begin{eqnarray*}
 \big| \theta - \theta^{\prime}  \big| & = & \big| ( I - \theta_{\mathcal S} \circ F_{t,r} \circ \theta^{-1}_{\mathcal S} ) (\theta) -  
 ( I - \theta_{\mathcal S} \circ F_{t,r} \circ \theta^{-1}_{\mathcal S} ) (\theta^{\prime}) \big|  \\
 & \leq &  C \scal{r}^{-\tau} \big| \theta - \theta^{\prime}  \big|   
\end{eqnarray*}
the second line following from  (\ref{operatornorm}) on the ball $ B (\theta_0, 2 \epsilon_{\omega_0}) $ which is convex. If $r$ is large enough, this implies that $ \theta = \theta^{\prime} $ hence that $ \omega = \omega^{\prime} $.

We next prove that $ F_{t,r} $ is surjective. More precisely, we show that if $r$ is large enough, then for all $ \omega \in \theta_{\mathcal S}^{-1} \big( B (\theta_0,\epsilon_{\omega_0}) \big) $ in the cover (\ref{recouvrement}), there exists $ \theta \in \theta_{\mathcal S}^{-1} \big( \overline{ B } (\theta_0,2 \epsilon_{\omega_0}) \big) $ such that
$$ \theta_{\mathcal S} (\omega) = \theta_{\mathcal S} \circ F_{t,r} \circ \theta_{\mathcal S}^{-1} (\theta) ,  $$
which we rewrite as the following fixed point equation
\begin{eqnarray}
 \theta =  T_{t,r} (\theta):= \big( I - \theta_{\mathcal S} \circ F_{t,r} \circ \theta_{\mathcal S}^{-1} \big) (\theta) + \theta_{\mathcal S} (\omega) . \label{pointfixePicard}
\end{eqnarray}
Indeed, we observe that the estimate (\ref{operatornorm}) still holds on $  \overline{ B } (\theta_0,2 \epsilon_{\omega_0}) $ by (\ref{conditionreconnaissance}) which implies that, for $r$ large enough, the map $ T_{t,r}  $ is $ 1/2 $-Lipschitz on $ \overline{B} (\theta_0,2 \epsilon_{\omega_0}) $. Furthermore, for $r$ large enough, (\ref{distanceabsolue}) implies that
$$ \big| \theta - \big( \theta_{\mathcal S} \circ F_{t,r} \circ \theta_{\mathcal S}^{-1} \big) (\theta) \big| \leq \epsilon_{\omega_0}, \qquad \theta \in  \overline{B} (\theta_0,2 \epsilon_{\omega_0}) , $$
hence that $ T_{t,r} $ maps $ \overline{B} (\theta_0,2 \epsilon_{\omega_0}) $ into $ \overline{B} (\theta_0,2 \epsilon_{\omega_0}) $, since $ |\theta_{\mathcal S}(\omega) - \theta_{0}| < \epsilon_{\omega_0} $. We can thus use the Picard fixed point Theorem to solve (\ref{pointfixePicard}) and this completes the proof of the surjectivity of $ F_{t,r} $.

All this shows that, for $r$ large enough, $ F_{t,r}  $ is (smooth and) bijective from $ {\mathcal S} $ to $ {\mathcal S} $. The smoothness of the inverse map follows from the inverse function theorem and (\ref{operatornorm}). More precisely, by (\ref{operatornorm}), we may assume for $r$ large enough that the differential of 
$ \theta_{\mathcal S} \circ F_{r,t} \circ \theta_{\mathcal S}^{-1} $ is invertible at any point of $ B (\theta_0, \epsilon_{\omega_0}) $ hence that $ \theta_{\mathcal S} \circ F_{r,t} \circ \theta_{\mathcal S}^{-1} $ is a local diffeomorphism close to any point of $ B (\theta_0, \epsilon_{\omega_0}) $. By (\ref{recouvrement}), we thus see that, for any $ \omega \in {\mathcal S} $, $ F_{t,r} $ is a diffeomorphism from a neighborhood of $ \omega $ onto a neighborhood  of $ F_{t,r} (\omega) $, which proves the smoothness of $ F_{t,r}^{-1} $. \finpreuve

{\it Author}: Jean-Marc Bouclet 

{\it Address}: Universit\'e  Toulouse 3,

$ $ Institut de Math\'ematique de Toulouse (UMR  CNRS 5219),

$ $ 118 route de Narbonne, F-31062 Toulouse Cedex 9

{\it Email}: Jean-Marc.Bouclet@math.univ-toulouse.fr


\begin{thebibliography}{99}

\bibitem{BouVP} {\sc J.-M. Bouclet}, {\it Absence of eigenvalue at the bottom of the continuous spectrum on asymptotically hyperbolic manifolds}, to appear in Annals of Global Analysis and Geometry.


\bibitem{Burq} {\sc N. Burq}, {\it Lower bounds for shape resonances widths of long range Schr\"odinger operators}, Amer. J. Math.
Vol. 124, Number 4 (2002), 677-735.

\bibitem{CaVo1} {\sc F. Cardoso, G. Vodev}, {\it High frequency resolvent estimates and energy decay of solutions to the wave equation},  Canad. Math. Bull.  47  (2004),  no. 4, 504-514.

\bibitem{CaVo2} {\sc \name}, {\it Uniform estimates of the resolvent of the Laplace-Beltrami operator on infinite volume Riemannian manifolds II},  Ann. Henri Poincar\'e  3  (2002),  no. 4, 673-691.

\bibitem{FaaDiBruno} {\sc G. M. Constantin, T. H. Savits}, {\it A multivariate Fa\`a Di Bruno formula with applications}, Trans. A.M.S., vol. 348, {\bf 2}, 503-520 (1996).

\bibitem{GHL} {\sc S. Gallot, D. Hulin, J. Lafontaine}, {\it Riemannian Geometry}, third edition, Springer (2004).




\bibitem{JoshiSaBaretto1} {\sc M. Joshi, A. S\'a Baretto}, {\it Recovering asymptotics of metrics from fixed energy scattering data},  Invent. Math.  137  (1999),  no. 1, 127-143.

\bibitem{JoshiSaBaretto2} {\sc \name}, {\it Inverse scattering on asymptotically hyperbolic manifolds},  Acta Math.  184  (2000),  no. 1, 41-86.

\bibitem{Kumura} {\sc H. Kumura}, {\it The radial curvature of an end that makes eigenvalues vanish in the essential spectrum I}, Math. Ann. (2010) 346:795-828.

\bibitem{Melrose} {\sc R. B. Melrose}, {\it Geometric Scattering Theory}, Cambridge Univ. Press (1995).


\end{thebibliography}
\end{document}